\documentclass{birkjour}
\usepackage{amsmath,amsthm}
\usepackage{amssymb, enumitem, mathrsfs, verbatim, bm}
\usepackage{color}

\catcode`\@=11 \@addtoreset{equation}{section}
\def\theequation{\thesection.\arabic{equation}}
\catcode`\@=12

\allowdisplaybreaks

\def\displ{\displaystyle}

\def\vec#1{\boldsymbol{#1}}

\def\bw{\boldsymbol{\omega}}

\def\bo{\vec{\omega}}
\def\S{\mathbb{S}}

\newcommand{\ep}{\varepsilon}
\newcommand{\vc}[1]{{\bf #1}}
\newcommand{\vu}{\vc{u}}

\newcommand{\vv}{\vc{v}}
\newcommand{\vz}{\vc{z}}
\newcommand\bph{\vec{\varphi}}

\newcommand{\Ds}{\mathbb{D}_x}
\newcommand{\Del}{\Delta_x}
\newcommand{\Dt}{\frac{\rm d}{{\rm d}t}}
\newcommand{\dt}{\,{\rm d} t }
\newcommand{\dx}{\,{\rm d} {x}}
\newcommand{\Div}{{\rm div}_x}
\newcommand{\Grad}{\nabla_x}
\newcommand{\intO}[1]{\int_{\Omega} #1 \ \dx}
\newcommand{\intTO}[1]{\int_0^T \int_\Omega #1 \ \dx\dt}

\def\O{\Omega}

\newcommand{\bfb}{\mathbf{b}}

\newcommand{\R}{\mathbb{R}}

\newtheorem{Theorem}{Theorem}[section]

\newtheorem{Corollary}[Theorem]{Corollary}

\theoremstyle{definition}
\newtheorem{Definition}[Theorem]{Definition}

\begin{document}

\title[Time-periodic weak solutions]{Time-periodic weak solutions to incompressible generalized Newtonian fluids}
\author[Anna Abbatiello]{Anna Abbatiello}
\address{Institute of Mathematics, Technische Universit\"{a}t Berlin, Stra{\ss}e des 17. Juni 136, 10623 Berlin, Germany}
\email{anna.abbatiello@tu-berlin.de}

\begin{abstract}
In this study we are interested in the Navier-Stokes-like system for generalized viscous fluids whose viscosity has a power-structure with exponent $q$. We develop an existence theory of periodic in time weak solutions to the three-dimensional flows subject to a periodic in time force datum whenever $q>\frac{6}{5}$, which is the optimal bound for the existence of weak solutions.  
\end{abstract}


\maketitle

\section{Introduction}
The Navier-Stokes-like system for an incompressible generalized viscous fluid reads 
\begin{subequations}\label{system-q}\begin{align}
\Div\,\vv=0 \mbox{ in } (0, T)\times\Omega,&\\
\partial_t\vv +\Div (\vv\otimes\vv) -\Div\, \S +\Grad p = \bfb \mbox{ in } (0, T)\times\Omega,&\label{eq-v}
\end{align}
where the viscous stress tensor $\S$ is a nonlinear function depending on the symmetric part of the velocity gradient 
$\mathbb{D}\vv :=\frac{1}{2}(\Grad\vv +(\Grad\vv)^t)$, and particularly we are concerned with the following rheological law:
\begin{equation}\label{eq-S}
\S=2\mu_0(\kappa+ |\mathbb{D}\vv|^2)^{\frac{q-2}{2}}\mathbb{D}\vv \mbox{ with } \mu_0>0, \kappa\geq 0, q\geq1,
\end{equation}
where for brevity we set $2\mu_0=1$. The Navier-Stokes system is the outcome by setting $q=2$ or equivalently by setting the generalized viscosity equal to a constant. 
The existence and the regularity theories to the Navier-Stokes-like system \eqref{system-q} give rise to the discussion on the bounds for the exponent or power-law index $q$. The best result for the existence of weak solutions to the three-dimensional initial and boundary value problem requires $q>{6}/{5}$ and it is obtained by employing the Lipschitz truncations method, see \cite{DRW}. 
However this bound can be removed in the class of \emph{generalized dissipative solutions}. They are 
 introduced  in \cite{AnnaEd} to a system describing the motion of a generalized viscous fluid
 and it is proved their long-time and large-data existence; particularly for a rheological law of type \eqref{eq-S} their existence holds for any $q\geq 1$. Despite the concept of generalized dissipative solution is weaker than the one of weak solution, it is showed in  \cite{AnnaEd} that these solutions enjoy the weak-strong uniqueness principle, in other words a \emph{generalized dissipative solution} must coincide with a classical solution as soon as and as long as the latter exists. 
Regarding the smoothness of the solutions to \eqref{system-q}, there are only partial answers. Especially the existence of smooth solutions corresponding to smooth data in the three-dimensional case is open. \\
In this paper we are concerned with the existence of periodic in time solutions to system \eqref{eq-v}-\eqref{eq-S} in a domain $\O\subset R^3$, completed with the boundary condition
\begin{align}\label{noslip}
\vv_{|\partial \Omega}=0.
\end{align}\end{subequations}
More specifically, we assume that the force datum is  a periodic in time function and we show that there exists a weak solution to system \eqref{system-q} which exhibits a reproductive property in time and has the same period of the force datum. We are interested in the degenerate case of the rheological law \eqref{eq-S} indeed we achieve the result with 
$$\kappa=0.$$
This question was first partially answered by Lions in \cite{Lions} where  it is assumed $q>3$ and $\S$ depends on the full gradient. In \cite{Crispo} the author considers the power-law exponent $q=q(x)$ a space-variable function but the convective term is neglected. In \cite{AbbatielloMaremonti} the difficulties are overcome proving the existence of \emph{regular solutions} under the assumption of smallness of the data, then it is established the existence of periodic in time solutions in a space-periodic domain with $q\in [\frac{5}{3}, 2)$. It is also relevant  to the problem the results in \cite{Barhoun} where the authors  requires $q\in (\frac{9}{5}, 2)$ and consider the space-periodic case. In this paper we prove the existence of periodic in time weak solutions for any $q$ in the whole range of values for which it is valid the existence of weak solutions.  
We construct a periodic in time approximations following the idea in \cite{Prouse63} used for the Navier-Stokes system. To this aim we  add a laplacian term to equation \eqref{eq-v}, then we perform a fixed point argument at the level of the Galerkin approximations and finally we take the limit in the approximations. The limit enjoys the reproductive property in time and is a weak solution to the system. We split the analysis into the two cases $q\geq 11/5$ and $q\in (6/5, 11/5)$. Once the result is achieved for $q\geq 11/5$, then we can construct time-periodic approximations in the case $q\in (6/5, 11/5)$ by adding a $p-$Laplacian term with $p=11/5$ to the system. (Note that a similar approximation scheme is exhibited in \cite{BleMalRaj}.) However, the case $q\in (6/5, 11/5)$ is more delicate and requires the introduction of one more approximation level in the generalized viscosity, in order to ensure the uniqueness of the solution to the Galerkin system which is needed for the fixed point argument. Finally we take the limit first in the Galerkin approximation, then in the parameter $\kappa$ in the generalized viscosity (achieving in this way the result for the degenerate case $\kappa=0$), and finally we get rid of the p-Laplacian and of the Laplacian regularizing terms through a trick which only use that we can keep a power of the parameter $\varepsilon$ in the estimates uniform in $\varepsilon$, as e.g. in \cite{BleMalRaj}.
In the end  the identification of the limit for the viscous stress tensor may be achieved through the standard Minty's trick when the limit of the velocity field itself can be used as test function in the weak formulation of \eqref{eq-v}, which is no longer possible when $q\in (6/5, 11/5)$. In this case we need to perform the Lipschitz truncations technique and more specifically the problem is solved employing the solenoidal version developed in \cite{BreDieSch} which simplifies this technique significantly. \\
The problem is here analyzed in the three-dimensional case for simplicity but the result may be achieved in
any spatial dimension $d$; then the supercritical case corresponds to $q \in (\frac{2d}{d+2}, 1+ \frac{2d}{d+2})$ and the
 subcritical/critical case to $q\geq 1+ \frac{2d}{d+2}$.

 \subsection{Time-periodic weak solution and main result}
\begin{Definition}
A field $\vc{f}:(0, +\infty)\times \O\to R^3$ is said to be {\emph{time-periodic}} with period $T>0$ if 
\begin{itemize}
\item[(i)]$\vc{f}\in L^\infty(0, T; X)$ with $(X, \|\cdot\|_X)$ a Banach space;
\item[(ii)]$\vc{f}(t+T, \cdot)=\vc{f}(t, \cdot)$ in $X$ for any $t\geq 0$.
\end{itemize}
\end{Definition}
\begin{Definition}
Let $T>0$, $q>6/5$ and let $\bfb\in C(0, +\infty; L^2(\O; R^3))$ a time-periodic function with period $T$. We say that $\vv$ is a {\emph{time-periodic weak solution}} to system \eqref{system-q} with period $T$ if the following hold:
\begin{itemize}
\item[(i)] $\vv\in L^\infty(0, T; L^2(\O; R^3))\cap L^q(0, T; W_{0, {\rm div}}^{1, q}(\O; R^3))$, \\
$\partial_t\vv\in (L^{q}(0,T; W_{0, {\rm div }}^{1, q}(\Omega; R^3)))^*$ if $q\geq 11/5$, \\
or $\partial_t\vv\in(L^{\frac{5q}{5q-6}}(0,T; W_{0, {\rm div }}^{1, \frac{5q}{5q-6}}(\Omega; R^3)))^*$ if $q\in (6/5, 11/5)$;
\item[(ii)]  there exists $\S\in L^{q'}((0, T)\times\O; R^{3\times 3})$ such that $$\S=|\Ds\vv|^{q-2}\Ds\vv \mbox{ a.e. in } (0, T)\times\O;$$
\item[(iii)] it holds the weak formulation
 \begin{equation*}\begin{split}
\int_0^T \langle\partial_t \vv, \bph\rangle \dt - \intTO{\vv\otimes\vv: \Grad \bph} + \intTO{\S:\Ds\bph}&\\
 = \intTO{\bfb \cdot\bph}
\mbox{ for any } \bph\in C^\infty((0, T); C_{0, {\rm div}}^\infty(\O; R^3)),&
 \end{split}\end{equation*}
 \item[(iv)] $\vv(T, \cdot)=\vv(0, \cdot)$ in $L^2(\O; R^3)$.
\end{itemize} 
Moreover, $\vv\in C(0, T; L^2(\O; R^3))$  if $q\geq 11/5$, or $\vv\in C_{\rm weak}(0, T; L^2(\O; R^3))$ if $q\in (6/5, 11/5)$.
\end{Definition}

We state the main result. 

 \begin{Theorem}\label{thm}
 Let $\O\subset R^3$ be a domain, let $T>0$ and let $q>6/5$. Assume $\bfb\in C(0, +\infty; L^2(\O; R^3))$ and that is a time-periodic function with period $T$. Then there exists a time-periodic weak solution $\vv$ with period $T$ to system \eqref{system-q}.
 \end{Theorem}
 
It is well known the property of extinction in a finite time of the solutions to a $q-$parabolic system when $q\in (1, 2)$, see \cite{CriGriMar}, \cite{DiB} and the same is true for a q-Stokes problem \cite{Abba}, \cite{Crispo}. In the following theorem we prove that in the case $q\in (6/5, 2)$ the system \eqref{system-q} admits a time-periodic weak solution which periodically becomes extinct. 
 
 \begin{Corollary}\label{corollary}
Let $\O\subset R^3$ be a domain, let $T>0$ and let $q\in (6/5, 2)$. Assume that $\bfb\in C(0, +\infty; L^2(\O; R^3))$ is a time-periodic function with period $T$ and that admits an extinction instant $\bar{t}\in (0, T)$, meaning that
$$\|\bfb(t)\|_{L^2(\Omega; R^3)}=0 \mbox{ for a.a. } t\in (\bar{t}, T),$$ 
and $\bar{t} + \frac{ \overline{K}^{2-q}}{\alpha (2-q)}\leq T$ where 
$$\overline{K}= \left( \frac{C_2(\max_{[0,T]}\|\bfb\|_{L^2(\Omega_T; R^3)}^{q'}+ 1)}{C_1C_S}\right)^{\frac{1}{q}}$$
and $C_1, C_2, C_S, \alpha$ are positive constant depending on $q$ and $\Omega$.
Then there exists a time-periodic weak solution $\vv$ with period $T$ to system \eqref{system-q} such that admits an extinction instant:
$$\|\vv(t)\|_{L^2(\Omega; R^3)}=0 \mbox{ for a.a. } t\in (\bar{t}_v, T)$$
with $\bar{t}_v\leq \bar{t} + \frac{ \overline{K}^{2-q}}{\alpha (2-q)}$.
\end{Corollary}
 
\paragraph{Notation and function spaces.} Let $\Omega\subset R^d$ a domain, i.e. an open bounded connected set. Throughout the paper we fix the dimension $d=3$, however we introduce the notation for a general dimension $d$. We denote the standard Lebesgue and Sobolev spaces of scalar functions by the usual notation $(L^r(\Omega), \|\cdot\|_{L^r(\Omega)})$ and $(W^{k, r}(\Omega),
 \|\cdot\|_{W^{k, r}(\Omega)})$. For vector-valued or tensor-valued functions, i.e. with value in $R^d$ or $R^{d\times d}$ respectively, we denote the corresponding Lebesgue and Sobolev spaces by  $(L^r(\Omega; R^d), \|\cdot\|_{L^r(\Omega; R^d)})$ and $(W^{k, r}(\Omega; R^{d\times d}), \|\cdot\|_{W^{k, r}(\Omega; R^{d\times d})})$ (with $1\leq r \leq +\infty$ and $k\in N$).  If $(X, \|\cdot\|_X)$ is a Banach space, $X^*$ is its dual space; then $C(0, T; X)$ is the relevant Bochner space and $C_{\rm weak}(0, T; X)$ is the space of weakly continuous function with value in $X$. Next, we define the space of compactly supported smooth functions and its subspace of solenoidal functions: 
\begin{align*}
C_{0}^\infty(\O; R^d)&:=\{ \vu:\Omega\to R^d, \vu \mbox{ smooth, } {\rm supp }\, \vu \subset \Omega\},\\
C_{0, {\rm div}}^\infty(\O; R^d)&:=\{ \vu\in C_0^\infty(\O; R^d), \Div \vu =0\},
\end{align*}
and their closures in $W^{k, r}$-norm for any $1<r<+\infty$ and $k\in N$: 
$$W_{0}^{k, r}(\O; R^d):=\overline{C_{0}^\infty(\O; R^d)}^{\|\cdot\|_{W^{k, r}}}, \ W_{0, {\rm div}}^{k, r}(\O; R^d):=\overline{C_{0, {\rm div}}^\infty(\O; R^d)}^{\|\cdot\|_{W^{k, r}}}.$$
Note that for a domain $\Omega$  (without further regularity assumption on the smoothness of $\partial \Omega$) we have the embedding $W_{0, {\rm div}}^{3, 2}(\O; R^3) \hookrightarrow  W_{0}^{3, 2}(\O; R^3) \hookrightarrow  W^{1, \infty}(\Omega; R^3)$. (We  will consider in particular the space $W_{0, {\rm div}}^{3, 2}(\O; R^3)$ for the basis for the Galerkin method.)
\\
As a consequence of the Poincar\'{e} and Korn inequalities the following norm are equivalent on the spaces $W_0^{1, r}$ and $W_{0, {\rm div}}^{1, r}$ for any $1<r<\infty$:
$$\|\Ds\bph\|_{L^r(\O; R^d)}\leq \|\Grad\bph\|_{L^r(\O; R^d)}\leq \|\bph\|_{W^{1, r}(\O; R^d)}\leq C_P C_K \|\Ds\bph\|_{L^r(\O; R^d)}$$
for all $\bph \in W^{1, r}(\O; R^d)$, where the constant $C_{P} >0$ that appears due to the Poincar\'e inequality depends on $r$ and $\Omega$, while the constant $C_{K} >0$  that appears due to the Korn inequality depends only on $r$.\\ Finally, for any vector valued function $\bph$ the symmetric part of the gradient is defined by $\Ds\bph:=\frac{1}{2}(\Grad\bph+(\Grad\bph)^t)$.

\section{Proof of Theorem \ref{thm}}

\subsection{The case $q\geq \frac{11}{5}$}

\paragraph{Existence of approximations.}
For any $\varepsilon>0$, let us introduce the approximating system 
\begin{subequations}\begin{align}
\partial_t\vv +\Div (\vv\otimes\vv) -\Div \S +\Grad p -\varepsilon \Del \vv = \bfb& \mbox{ in } (0, T)\times\Omega,\\
\Div\,\vv=0& \mbox{ in } (0, T)\times\Omega,\\
\vv=0& \mbox{ on } (0, T)\times\partial\Omega.
\end{align}\end{subequations}
where  $$\S=|\Ds\vv|^{q-2}\Ds\vv.$$
In order to fix the basis for the Galerkin method, consider the eigenvalue problem
\begin{equation}
(\!(\bw,\bph)\!)=\lambda (\bw, \bph) \ \mbox{ for all } \bph\in W^{3,2}_{0, {\rm div}}(\Omega; R^3)
\end{equation}
with $\lambda\in R$ and $\bw\in W^{3,2}_{0, {\rm div}}(\Omega; R^3)$, and where $(\cdot,\cdot)$ is the scalar product in $L^2(\Omega; R^3)$ while $(\!(\cdot,\cdot)\!)$ is the scalar product in $W^{3,2}_{0, {\rm div}}(\Omega; R^3)$ defined as $(\!(\bw,\bph)\!):=(\nabla^3\bw,\nabla^3\bph)+(\bw,\bph)$. It is known (see e.g. \cite[Appendix A.4]{Malekbook}) that there exist eigenvalues $\{\lambda_m\}_{m\in N}$ and corresponding eigenfunctions $\{ \bw^m\}_{m\in N}$  such that they are orthonormal in $L^2(\Omega, R^3)$ and orthogonal in $W^{3,2}_{0, {\rm div}}$. For any fixed $n\in N$ the projectors $P^n:W^{3,2}_{0, {\rm div}}\to X_n:=\mbox{span}\{\bw^1, \dots, \bw^n\}$ defined as $P^n\vv:=\sum_{i=1}^n (\vv, \bw^i)\bw^i$ are continuous orthonormal projectors in $L^2(\Omega, R^3)$.  It is worth recalling that $W^{3,2}_{0, {\rm div}}(\Omega; R^3)\hookrightarrow W^{3,2}_{0}(\Omega; R^3)\hookrightarrow W^{1,\infty}(\Omega, R^3).$ \\
Now, given $\vv_0^n\in \mbox{span}\{\bw^1, \dots, \bw^n\}$ with
\begin{equation}\label{assumption0}
\|\vv_0^n\|_{L^2(\Omega; R^3)}\leq K
\end{equation}
with $K$ constant defined as
$$ K:= \left( \frac{C_2\max_{[0,T]}\|\bfb\|_{L^2(\O; R^3)}^{q'}}{C_1C_S}\right)^{\frac{1}{q}}$$ 
(with constants $C_1, C_2, C_S$ defined later in the proof, also the reasons will be understood), we look for the Galerkin approximation 
$$\vv^n(t,x):=\sum_{i=1}^n c^n_i(t)\bw^i(x)$$
 such that it satisfies for any $k=1,\dots, n$
\begin{subequations}\begin{equation}
\!(\partial_t\vv^n, \bo^k)- (\vv^n\otimes \vv^n, \Grad \bw^k) +  (\S^n, \Ds\bw^k) 
+\varepsilon(\Grad \vv^n, \Grad\bw^k) =  (\bfb, \bw^k) \label{system0}\end{equation}
where
\begin{equation}
\S^n:=  |\Ds\vv^n|^{{q-2}}\Ds\vv^n,
\end{equation}
and 
\begin{equation}
\vv^n(0)=\vv_0^n,
\end{equation}
 \end{subequations}
 meaning in components that for any $k=1,\dots, n$ 
\begin{subequations}\label{ode0}\begin{align}
 (c^n_k(t))'&= \!\sum_{i,j=1}^n c^n_ic^n_j f_{ijk}-\left(\S^n, \Ds\bw^k\right) - \varepsilon\! \sum_{i=1}^n c^n_i g_{ik}+b_k, 
 \label{ode1-0}\\
c^n_k(0)&=c^n_{0, k}:=(\vv_0^n, \bw^k),
\end{align}\end{subequations}
where we used the notation 
\begin{equation*}
f_{ijk}:=(\bw^i\otimes \bw^j, \Grad \bw^k),\
 g_{ik}:= (\Grad \bw^i, \Grad\bw^k), \ 
 b_k:=(\bfb, \bw^k).
 \end{equation*}
 The existence and uniqueness of the solution ${\bf c}^n(t)=(c^n_1(t), \dots, c^n_n(t))$ to the Cauchy problem \eqref{ode0} in a local time interval $[0, t_n)$ follow by standard results on ODEs being the right-hand side of \eqref{ode1-0} a Lipschitz function thanks to the assumption $q\geq{11}/{5}$. \\
Multiplying \eqref{system0} by $c^n_k$ and summing over $k=1,\dots, n$ we get
\begin{equation}\label{diff-en-0}\begin{array}{l}\displaystyle\vspace{6pt}
\frac{d}{dt} \|\vv^n(t)\|_{L^2(\O; R^3)}^2+ C_1\|\Ds\vv^n\|_{L^q(\Omega; R^{3\times 3})}^q+\varepsilon\|\Grad\vv^n\|_{L^2( \Omega; R^{3\times 3})}^2 \\\displaystyle
\hfill\leq C_2\|\bfb\|_{L^2(\O; R^3)}^{q'},\end{array}
\end{equation}
then integrating in time over $(0,t)$ we obtain
\begin{equation}\label{diff-en2-0}
\!\!\|\vv^n(t)\|_{L^2(\O; R^3)}^2\!\!\leq \!\|\vv_0^n\|_{L^2(\O; R^3)}^2+ C_2\,T\,\max_{[0,T]}\|\bfb\|_{L^2(\O; R^3)}^{q'}  \mbox{ for all } t\!\in\![0,T].
\end{equation}
Note that $\|\vv^n(t)\|_{L^2(\O; R^3)}^2=|{\bf c}^n(t)|^2$,  thus  $\vv^n(t)$ is well-defined on the whole interval $[0,T]$. 
\\

\paragraph{Existence of a periodic in time approximating solution.}
Now, by \eqref{diff-en-0} and by the assumption in \eqref{assumption0} we realize that 
\begin{equation}\label{R}
|{\bf c}^n(t)|=\|\vv^n(t)\|_{L^2(\O; R^3)}\leq K \mbox{ for all } t\in [0,T],
\end{equation}
indeed if there exists an instant $\bar{t}$ such that $\|\vv^n(\bar{t})\|_2= K$
then 
\begin{equation}\label{max}
\frac{d}{dt} \|\vv^n(\bar{t})\|_{L^2(\O; R^3)}^2< C_2\max_{[0,T]}\|\bfb\|_{L^2(\O; R^3)}^{q'}- C_1C_SK^{q}=0
\end{equation}
where $C_S$ is the constant due to the embedding $W_0^{1, q} \hookrightarrow L^2$ and to the Korn inequality. 
Being $\vv^n\in C(0, T; L^2(\Omega; R^3))$,  \eqref{max} implies that $\|\vv^n(t)\|_{L^2(\O; R^3)}^2$ is a non-increasing function whenever $\|\vv^n(t)\|_{L^2(\O; R^3)}= K$, thus \eqref{R} holds. 

Consider ${\bf c}^n(t)=(c^n_1(t), \dots, c^n_n(t))$ and ${\bf d}^n(t)=(d^n_1(t), \dots, d^n_n(t))$ solutions of \eqref{ode1-0} corresponding to the same force $\bfb=\bfb(t, x)$ and set $\alpha^n_k(t):= c^n_k(t)- d^n_k(t)$ for any $k=1,\dots, n$. It holds for any $k=1,\dots,n$
\begin{equation}
(\alpha^n_k(t))' - \!\sum_{i,j,k=1}^n (\alpha^n_ic^n_j + d^n_i\alpha^n_j)f_{ijk}+\left(\S^{n}_1-\S^{n}_2, \Ds\bw^k \right) \\ \displaystyle\vspace{6pt}
+\varepsilon \sum_{i=1}^n\alpha^n_i g_{ik} =0 
\end{equation}
where $$\S^{n}_1:=\sum_{i=1}^n c^n_i \,\left|\sum_{j=1}^nc^n_j\Ds\bw^j\right|^{q-2}\!\!\!\!\Ds\bw^i \ \mbox{ and } \ \S^{n}_2:=\sum_{i=1}^n d^n_i\, \left |\sum_{j=1}^nd^n_j\Ds\bw^j\right|^{q-2}\!\!\!\!\Ds\bw^i.$$
Multiplying by $\alpha^n_k$, summing over $k=1,\dots,n$ and employing the monotonicity of the operator in \eqref{eq-S} and the Poincar\'e inequality, we get
\begin{equation}\label{diff3-0}\begin{array}{l}\displaystyle\vspace{6pt}
\frac{1}{2}\frac{d}{dt} \|\vu^n(t)\|_{L^2(\O; R^3)}^2 - \sum_{i,j,k=1}^n \alpha^n_i(t)c^n_j(t) \alpha^n_k(t)f_{ijk} \\ \displaystyle
\hfill +\varepsilon\, C_3\|\vu^n(t)\|_{L^2(\O; R^3)}^2\leq 0
\end{array}\end{equation}
where 
$$\vu^n(t):=\sum_{k=1}^n \alpha^n_k(t)\bw^k.$$
Now, we follow the idea in \cite{Prouse63} to treat the convective term.
Employing the regularity of the basis functions $\bw^i$, for any $n\in N$ it holds that 
\begin{equation}\label{convective0}\begin{split}
&\!\!\!\left|\sum_{i,j,k=1}^n \!\!\!\alpha^n_i(t)c^n_j(t) \alpha^n_k(t)f_{ijk}\right| \leq n \max_{i,j,k} f_{ijk} \max_{[0,T]} |{\bf c}^n(t)|\sum_{i,k=1}^n |\alpha^n_i(t)| |\alpha^n_k(t)|\\
&\leq n^2 \max_{i,j,k} f_{ijk} \max_{[0,T]} |{\bf c}^n(t)| \left(\sum_{i=1}^n\frac{1}{2} |\alpha^n_i(t)|^2+ \sum_{k=1}^n\frac{1}{2} |\alpha^n_k(t)|^2\right)\\
&=C_4 \|\vu^n(t)\|_{L^2(\O; R^3)}^2
\end{split}\end{equation}
where we set $C_4:=n^2 \max_{i,j,k} f_{ijk} \max_{[0,T]} |{\bf c}^n(t)|$.
Employing \eqref{convective0} in \eqref{diff3-0}, we obtain 
\begin{equation}
\frac{d}{dt}\|\vu^n(t)\|_{L^2(\O; R^3)}\leq (C_4-\varepsilon\, C_3) \|\vu^n(t)\|_{L^2(\O; R^3)},
\end{equation}
which implies 
\begin{equation}
\|\vu^n(t)\|_{L^2(\O; R^3)}\leq e^{(C_4-\varepsilon\, C_3)t} \|\vu^n(0)\|_{L^2(\O; R^3)} \mbox{ for any } t\in [0,T],
\end{equation}
and it means that 
\begin{equation}\label{continuity0}
\|\vv^n(t)-\vz^n(t)\|_{L^2(\O; R^3)}\leq e^{(C_4-\varepsilon\, C_3)t} \|\vv^n(0)-\vz^n(0)\|_{L^2(\O; R^3)} 
\end{equation}
for any $t\in [0,T]$, where $\vv^n(t):=\sum_{k=1}^n c^n_k(t)\bw^k$ and $\vz^n(t):=\sum_{k=1}^n d^n_k(t)\bw^k.$
Let us consider the mapping
$$ F: R^n\to R^n, \, F({\bf c}^n_0):={\bf c}^n(T)$$
where ${\bf c}^n(t)$ is the solution to the system \eqref{ode0} corresponding to the initial value ${\bf c}^n_0$. The mapping $F$ is well-defined as the solution to \eqref{ode0} there exists and it is unique. Moreover, $F$ transforms the ball $B_{K}$ of $R^n$ into itself by \eqref{R} and it is continuous by \eqref{continuity0}. Then the Brouwer theorem ensures the existence of a fixed point, hence for any $n\in N$ there exists $\vv^n(t)$ solution to \eqref{system0} such that 
$$\vv^n(T)=\vv^n(0).$$

\paragraph{ Limit as $n\to +\infty$.} 
Integrating \eqref{diff-en-0} in time over $(0,T)$ we have 
\begin{equation}\label{uniform1-0}
\sup_n\|\Ds\vv^n\|_{L^q(\Omega_T; R^{3\times 3})}+\varepsilon\|\Grad\vv^n\|_{L^2(\Omega_T; R^{3\times 3})}<+\infty
\end{equation}
where we introduced the notation $\O_T:= (0, T)\times \O$, and as a consequence 
\begin{equation}
\sup_n\|\S^n\|_{L^{q'}(\Omega_T; R^{3\times 3})}<+\infty.
\end{equation}
From the interpolation inequality 
$$ \intO{|\vv^n|^{\frac{5q}{3}}}\leq \left(\intO{|\vv^n|^2}\right)^{\frac{q}{3}}\left(\intO{|\vv^n|^{\frac{3q}{3-q}}}\right)^{\frac{3-q}{3}},$$
integrating in time and using Sobolev embeddings, we achieve that
\begin{equation}\label{5q-3}
\intTO{|\vv^n|^{\frac{5q}{3}}} \leq \sup_{[0,T]} \|\vv^n\|_{L^2(\Omega; R^3)}^{\frac{2q}{3}}\intTO{|\Grad\vv^n|^q},
\end{equation}
next \eqref{uniform1-0} and \eqref{diff-en2-0} together with the Korn inequality imply that
\begin{equation}\label{uniform2-0}
 \sup_n\|\vv^n\|_{L^{\frac{5q}{3}}(\Omega_T; R^3)}<+\infty.
\end{equation}
For any $\bph\in L^q(0,T; W_{0, {\rm div}}^{3, 2}(\Omega; R^3))$ 
$$\int_0^T\int_\Omega \partial_t\vv^n\cdot \bph \dx\dt=\int_0^T\int_\Omega \partial_t\vv^n\cdot P^n\bph\dx\dt$$
then from \eqref{system0} and by virtue of \eqref{diff-en2-0}, \eqref{uniform1-0} and \eqref{uniform2-0}, we obtain that
\begin{equation}
\sup_n\|\partial_t\vv^n\|_{(L^{q}(0,T; W_{0, {\rm div}}^{3, 2}(\Omega ; R^3)))^*}<+\infty. 
\end{equation}
For any fixed $\varepsilon>0$ and for any $n\in N$ the established estimates uniform respect to $n$ yield the existence of $\vv$ and $\S$  such that the following convergences hold (for suitable subsequences not relabelled):
\begin{align}
\vv^n &\to \vv &&\mbox{ weakly-* in } L^\infty(0,T; L^2(\Omega; R^3)),\label{convergence10}\\
\vv^n &\to \vv &&\mbox{ weakly  in } L^q(0,T; W^{1, q}_{0, {\rm div}}(\Omega; R^3)), \label{debole}\\
\partial_t \vv^n &\to \partial_t\vv &&\mbox{ weakly in } (L^{q}(0,T; W_{0, {\rm div}}^{3, 2}(\Omega; R^3)))^*, \label{time}\\
\S^n &\to \S &&\mbox{ weakly  in } L^{q'}(\Omega_T; R^{3\times 3}),\label{convergence50}\\
\Grad \vv^n &\to \Grad\vv &&\mbox{ weakly  in } L^2(\Omega_T; R^{3\times 3}). \label{conv-grad}
\end{align}
Furthermore, using \eqref{uniform2-0} and \eqref{time}, the Aubin-Lions lemma applied to $W_{0, {\rm div}}^{3, 2} \hookrightarrow W_0^{1, q}\hookrightarrow\hookrightarrow L^q \hookrightarrow (W_{0, {\rm div}}^{3, 2})^*$ gives
$$\vv^n\to \vv \mbox{ strongly in } L^q(\Omega_T; R^3),$$
which together with the following interpolation inequality
\begin{equation}\label{aubin2}
 \|\vv^n-\vv\|_r\leq \|\vv^n-\vv\|_q^a\|\vv^n-\vv\|_{\frac{5q}{3}}^{1-a}\leq C\|\vv^n-\vv\|_q^a \mbox{ for any } r< \frac{5q}{3}
 \end{equation}
implies  
\begin{align}\label{strong0}
\vv^n \to \vv \mbox{ strongly in } L^r(\Omega_T; R^3) \mbox{ for any } r< \frac{5q}{3}.
\end{align}
It is worth noting that \eqref{convergence10}, \eqref{time} ensure
\begin{equation}
\vv\in C(0,T; L^2(\Omega; R^3)).
\end{equation}
Moreover, the Ascoli-Arzel\`{a} theorem implies that for any $\bph\in C^\infty_{0, {\rm div}}(\Omega; R^3)$
\begin{equation}\label{uniform-weak0}
(\vv^n(t), \bph)\to (\vv(t), \bph) \mbox{ uniformly in } t\in[0,T],  \mbox{ as } n\to+\infty.
\end{equation}
On the other side, by \eqref{assumption0} we get the existence of a subsequence $\vv_0^n$ such that 
\begin{equation}\label{weak0}
\vv_0^n\to \vv_0 \mbox{ weakly  in } L^2(\Omega; R^3)
\end{equation}
with  $\vv_0$ satisfying again the bound \eqref{assumption0}.
Therefore, employing \eqref{uniform-weak0} and \eqref{weak0} it follows
\begin{equation}
(\vv(T)-\vv_0, \bph)=\lim_{n\to+\infty} (\vv^n(T)-\vv_0^n, \bph)=0
\end{equation}
i.e. the limit function $\vv$ is time-periodic with period $T$. Finally thanks to the convergences \eqref{convergence10}--\eqref{strong0},  using that $W_{0, {\rm div}}^{3, 2}$ is dense in $W_{0, {\rm div}}^{1, q}$ for any $q\geq 1$, and employing  a standard density argument  we achieve that $\vv$ and $\S$ enjoy the following weak formulation: 
\begin{equation}\label{system2}
\begin{array}{l}\displaystyle\vspace{6pt}
\int_0^T\langle \partial_t\vv, \bph\rangle\dt+\int_0^T \int_\Omega \S :\Ds\bph\dx\dt+\varepsilon\int_0^T\int_\Omega\Grad\vv\!:\!\Grad\bph\dx\dt
\\ \displaystyle\vspace{6pt}\hfill 
= \int_0^T\int_\Omega \vv\otimes\vv:\Grad\bph\dx\dt +\intTO{\bfb\cdot \bph} \end{array}
\end{equation}
for any  $\bph\in L^{q}(0,T; W_{0, {\rm div}}^{1, q}(\Omega; R^3))$. Thus $\partial_t\vv\in (L^{q}(0,T; W_{0, {\rm div}}^{1, q}(\Omega; R^3)))^*$. 
It remains to show that $\S$ fulfils \eqref{eq-S}. Multiplying \eqref{system0} by $c^n_k$, summing over $k=1,\dots, n$ and integrating in time over $(0, T)$ we have
 $$ \int_0^T\int_\Omega \S^n :\Ds\vv^n\dx\dt +\varepsilon\! \int_0^T\int_\Omega |\Grad\vv^n|^2\dx\dt=\intTO{\bfb\cdot\vv^n},$$
 taking the limsup and employing  \eqref{debole} and \eqref{conv-grad}  it follows
  \begin{equation}\label{limsup}
   \limsup_{n\to +\infty} \int_0^T\int_\Omega \S^n:\Ds\vv^n\dx\dt\leq -\varepsilon\|\Grad\vv\|_2^2+\intTO{\bfb\cdot\vv}.
  \end{equation}
Comparing \eqref{limsup} with the outcome of \eqref{system2} choosing $\bph=\vv$ it follows
 \begin{equation}\label{limsup0}
 \limsup_{n\to +\infty} \int_0^T\int_\Omega \S^n:\Ds\vv^n\dx\dt\leq  \int_0^T\int_\Omega \S :\Ds\vv\dx\dt
 \end{equation}
 and finally by virtue of the established convergences together with \eqref{limsup0} the standard Minty's trick ensures that 
 $$\S= |\mathbb{D}\vv|^{q-2}\mathbb{D}\vv \mbox{ a.e. in } \Omega_T.$$

\paragraph{Limit as $\varepsilon\to 0$.} For any $\varepsilon>0$ we constructed $\vv^\varepsilon$  and $\S^\varepsilon$ fulfilling 
\begin{equation}\label{system3}
\begin{array}{l}\displaystyle\vspace{6pt}
\!\!\!\!\!\!\!\int_0^T\langle \partial_t\vv^\ep, \bph\rangle\dt+\int_0^T \int_\Omega \S^\ep :\Ds\bph\dx\dt+\varepsilon\int_0^T\int_\Omega\Grad\vv^\ep:\Grad\bph\dx\dt
\\ \displaystyle\vspace{6pt}
= \int_0^T\int_\Omega \vv^\ep\otimes\vv^\ep:\Grad\bph\dx\dt +\intTO{\bfb\cdot \bph} \end{array}
\end{equation}
for any $\bph\in L^{q}(0,T; W_{0, {\rm div}}^{1, q}(\Omega; R^3))$, with $\vv^\ep$ time-periodic with period $T$. Taking $\vv^\varepsilon$ as test function we derive uniform estimates that ensure the following convergences as $\ep\to 0$ (for suitable subsequences not relabelled): 
 \begin{align}
\vv^\varepsilon &\to \vv &&\mbox{ weakly-* in } L^\infty(0,T; L^2(\Omega; R^3)),\label{convergence100}\\
 \vv^\varepsilon &\to \vv &&\mbox{ weakly  in } L^q(0,T; W^{1, q}_0(\Omega; R^3)), \label{debole0}\\
\partial_t\vv^\varepsilon &\to \partial_t\vv &&\mbox{ weakly in } (L^{q}(0,T; W_0^{1, q}(\Omega; R^3)))^*, \label{time0}\\
\vv^\varepsilon&\to \vv &&\mbox{ strongly in } L^r(\Omega_T; R^{3}) \mbox{ for any } r< \frac{5q}{3},\label{strong00}\\
\S^\varepsilon &\to \S &&\mbox{ weakly  in } L^{q'}(\Omega_T; R^{3\times 3}).\label{convergence500}
\end{align}
Being $5q/6 < q$, $\|\Grad\vv^\varepsilon\|_{\frac{5q}{6}}$ is uniformly bounded and as a consequence it results
\begin{equation}\label{epsilon}
\left|\varepsilon\!\!\int_0^T\!\!\!\int_\Omega\nabla\vv\!:\!\nabla\bph \right|\leq \varepsilon \|\nabla\vv^\varepsilon\|_{\frac{5q}{6}}\|\nabla\bph\|_{\frac{5q}{5q-6}}\to 0 \mbox{ as } \varepsilon \to 0.
\end{equation}
We are in a position to take the limit as $\ep\to 0$ in \eqref{system3} and employing the convergences in \eqref{convergence100}-\eqref{epsilon} we obtain the existence of $(\vv, \S)$ such that 
\begin{equation}\label{system-w}
\begin{array}{l}\displaystyle\vspace{6pt}
\int_0^T\langle \partial_t\vv, \bph\rangle\dt+\int_0^T \int_\Omega \S :\Ds\bph\dx\dt
= \int_0^T\int_\Omega \vv\otimes\vv:\Grad\bph\dx\dt
\\ \displaystyle\vspace{6pt}\hfill  +\intTO{\bfb\cdot \bph} \mbox{  for any } \bph\in L^{q}(0,T; W_{0, {\rm div}}^{1, q}(\Omega; R^3)).
\end{array}
\end{equation}
Employing the same tools of the previous passage to the limit as $n\to + \infty$ we obtain that $\vv$ is time-periodic with period $T$ and that $\S$ fulfils the rheological law \eqref{eq-S}.   

\subsection{The case $q\in (\frac{6}{5}, \frac{11}{5})$}

\paragraph{Existence of a periodic in time approximating solution.}
For any $\ep, \kappa >0$ let us introduce the approximating system 
\begin{subequations}\label{system-k-ep}\begin{align}
\!\!\!\partial_t \vv+\Div (\vv\otimes\vv) -\Div \S +\Grad p -\varepsilon (\Del \vv+ \Div (|\Ds\vv|^{\frac{1}{5}} \Ds\vv))&= \bfb 
\\
\Div\,\vv&=0
\end{align}
with 
\begin{equation}
\S=(\kappa+ |\Ds\vv|^2)^{\frac{q-2}{2}}\Ds\vv
\end{equation}
 and $(t, x)\in(0, T) \times \Omega$ and boundary condition
 \begin{equation}
 \vv=0 \mbox{ on } (0, T) \times \partial\Omega.
 \end{equation}
\end{subequations}
It is worth remarking that in this section, being $q\in (\frac{6}{5}, \frac{11}{5})$, we need to approximate the degenerate stress tensor with a non-degenerate one in order to ensure the uniqueness of the Galerkin approximation, which is crucial for our argument. (Subsequently we will let $k\to 0$.) Moreover once we added the $p$-Laplacian term with $p=11/5$ we may conclude the existence of time-periodic Galerkin approximations to the approximating system \eqref{system-k-ep} by virtue of the results proved in the former section; however for the sake of clarity and for the reader's convenience we perform the main steps below.\\
Given $\vv_0^n\in \mbox{span}\{\bw^1, \dots, \bw^n\}$ such that
\begin{equation}\label{assumption1}
\|\vv_0^n\|_{L^2(\Omega; R^3)}\leq K:= \left( \frac{C_2(\max_{[0,T]}\|\bfb\|_{L^2(\Omega_T; R^3)}^{q'}+\kappa^{\frac{q}{2}})}{C_1C_S}\right)^{\frac{1}{q}}
\end{equation}
with constants $C_1, C_2, C_S$ defined later in the proof, we look for the Galerkin approximation 
$$\vv^n(t,x):=\sum_{i=1}^n c^n_i(t)\bw^i(x)$$ 
such that it satisfies for any $k=1,\dots, n$
\begin{subequations}\begin{align} 
&\begin{array}{l}\displaystyle\vspace{6pt}
(\partial_t \vv^n, \bo^k)+  (\S^n - \vv^n\otimes \vv^n, \Ds\bw^k) +\varepsilon(\Grad \vv^n, \Grad\bw^k)\\
\displaystyle\vspace{6pt} \hfill+\varepsilon (|\Ds\vv^n|^{\frac{1}{5}}\Ds\vv^n, \Ds\bw^k)=  (\bfb, \bw^k), \end{array}
\label{system}\\
&\vv^n(0)=\vv_0^n,
\end{align}
where 
\begin{equation}
\S^n:=(\kappa +|\Ds\vv^n|^2)^{\frac{q-2}{2}}\Ds\vv^n, 
\end{equation}
\end{subequations}
 meaning in components that for any $k=1,\dots, n$
\begin{subequations}\label{ode}\begin{align}
(c^n_k(t))'&= \sum_{i,j=1}^n c^n_ic^n_j f_{ijk}-(\S^n, \Ds\bw_k) - \varepsilon \sum_{i=1}^n c^n_i g_{ik}\label{ode1}\\
&-\varepsilon (|\Ds\vv^n|^{\frac{1}{5}} \Ds\vv^n, \Ds\bw^k) + b_k, \nonumber\\
c^n_k(0)&=(\vv_0^n, \bw^k),
\end{align}\end{subequations}
and where we introduced the notation 
$$f_{ijk}:=(\bw^i\otimes \bw^j, \Grad \bw^k), \ g_{ik}:= (\Grad \bw^i, \Grad\bw^k),  \ b_k:=(\bfb, \bw^k).$$
The existence of a solution ${\bf c}^n(t):= (c^n_1, \dots, c^n_n)$ to the Cauchy problem \eqref{ode} in a local time interval $[0, t_n)$ follows by standard results on ODEs. Thanks to the presence of $\kappa>0$, we can state that the right-hand side of \eqref{ode1} is a Lipschitz function and thus the uniqueness of ${\bf c}^n(t)$ follows. 
Now, let us multiply \eqref{system} by $c^n_k$ and sum over $k=1,\dots, n$ 
\begin{equation}\label{diff-en}\begin{array}{l}\displ\vspace{8pt}
\Dt\|\vv^n(t)\|_{L^2(\O; R^3)}^2+ C_1\|\Ds\vv^n\|_{L^q(\Omega; R^{3\times 3})}^q +\varepsilon\|\Grad\vv^n\|_{L^2(\Omega; R^{3\times 3})}^2\\ \displ\vspace{8pt}\hfill
+\varepsilon\|\Ds\vv^n\|_{L^{\frac{11}{5}}(\Omega; R^{3\times 3})}^{\frac{11}{5}} \leq C_2\|\bfb\|_{L^2(\O; R^3)}^{q'}+C_2\kappa^{\frac{q}{2}},
\end{array}\end{equation}
then let us integrate in time over $(0,t)$, we obtain for all  $t\in[0,T]$
\begin{equation}\label{diff-en2}
\|\vv^n(t)\|_{L^2(\O; R^3)}^2\leq \|\vv_0^n\|_{L^2(\O; R^3)}^2+ C_2T\max_{[0,T]}\|\bfb\|_{L^2(\O; R^3)}^{q'}+C_2\kappa^{\frac{q}{2}}T, 
\end{equation}
where $\|\vv^n(t)\|_{L^2(\O; R^3)}^2=|{\bf c}^n(t)|^2$. This implies that $\vv^n(t)$ is well-defined on the whole interval $[0,T]$. 
\\
Let us introduce the mapping
$$ F:\R^n\to\R^n, \, F({\bf c}^n_0):={\bf c}^n(T)$$
where ${\bf c}^n(t)$ is the unique solution to the system \eqref{ode} corresponding to the initial value ${\bf c}^n_0$. Repeating the arguments showed in the analogous paragraph in the case $q\geq 11/5$, we obtain that $F$ transforms the ball $B_{K}$ of $R^n$ of radius $K$ into itself and it is continuous. Therefore by the Brouwer theorem it follows the existence of a fixed point, which means that for any $n\in N$ there exists $\vv^n(t)$ solution to \eqref{system} such that $\vv^n(T)=\vv^n(0).$
\\

\paragraph{\bf Limit as $n\to + \infty$.} 
Integrating \eqref{diff-en} in time over $(0,T)$ we have 
\begin{equation}\label{uniform1}\begin{split}
\sup_{ n\in N}&\left(\|\Ds\vv^n\|_{L^q(\Omega_T; R^{3\times 3})}+\|\varepsilon^{\frac{1}{2}}\Grad\vv^n\|_{L^2(\Omega_T; R^{3\times 3})}\right.
\\&+\left. \|\varepsilon^\frac{5}{11}\Ds\vv^n\|_{L^{\frac{11}{5}}(\Omega_T; R^{3\times 3})}\right)<+\infty,
\end{split}\end{equation}
and as a consequence 
\begin{equation}\label{uniform-S}
\sup_{ n\in N}\|\S^n\|_{L^{q'}(\Omega_T;  R^{3\times 3})}<+\infty.
\end{equation}
Analogously to \eqref{5q-3}, it holds  
$$
\int_0^T \int_\Omega |\vv^n|^{\frac{11}{3}}\dx\dt \leq \sup_{[0,T]}\left( \int_\O|\vv^n|^2\dx\right)^{\frac{11}{15}}\int_0^T\int_\O |\Grad\vv^n|^{\frac{11}{5}}\dx\dt,
$$
then the Korn inequality, \eqref{diff-en2} and \eqref{uniform1}  imply that
\begin{equation}\label{uniform2}
\sup_{ n\in N} \|\vv^n\|_{L^{\frac{11}{3}}(\Omega_T; R^3)}<+\infty.
\end{equation}
By virtue of estimates \eqref{diff-en2}, \eqref{uniform1} and \eqref{uniform2} we derive that the time derivative of $\vv^n$ is uniformly bounded in the following space:
\begin{equation}\label{uniform-t}
\sup_{ n\in N}\|\partial_t\vv^n\|_{(L^{\frac{11}{5}}(0,T; W_0^{1, \frac{11}{5}}(\Omega; R^3)))^*}<+\infty. 
\end{equation}
The obtained uniform estimates yield the existence of $\vv$ and $\S$ such that, for suitable subsequence not relabelled, the following convergences hold:
\begin{align}
&\vv^n\to \vv \mbox{ weakly-* in } L^\infty(0,T; L^2(\Omega; R^3)),\label{convergence1}\\
& \vv^n\to \vv \mbox{ weakly  in } L^\frac{11}{5}(0,T; W^{1, \frac{11}{5}}_{0, {\rm div}}(\Omega; R^3)),\\
&\partial_t\vv^n\to \partial_t\vv \mbox{ weakly-* in } (L^{\frac{11}{5}}(0,T; W_{0, {\rm div}}^{1, \frac{11}{5}}(\Omega; R^3)))^*,\\
&\S^n\to \S\mbox{ weakly  in } L^{q'}(\Omega_T; R^{3\times 3}).\label{convergence5}
\end{align}
The Aubin-Lions compactness lemma applied to $W_0^{1, \frac{11}{5}}\hookrightarrow\hookrightarrow L^\frac{11}{5} \hookrightarrow (W_0^{1, \frac{11}{5}})^*$ and the interpolation inequality in \eqref{aubin2} with $q=11/5$ imply that
\begin{equation}
\vv^n\to \vv \mbox{ strongly in } L^r(\Omega_T; R^3) \mbox{ for any } r< \frac{11}{3}.\label{strong}
\end{equation}
On the other side, from \eqref{assumption1} it follows that for a suitable subsequence it holds
\begin{equation}
\vv^n_0\to \vv_0 \mbox{ weakly in } L^2(\O; R^3) \mbox{ and } \|\vv_0\|_{L^2(\O; R^3)}\leq K.
\end{equation}
Now, we can follow step by step the  passage to the limit as $n\to +\infty$ in the section $q\geq 11/5$, and we obtain $\vv$ time-periodic weak solution to system \eqref{system-k-ep} with period $T$. 
\\

\paragraph{Limit as $\kappa\to0$.}
Therefore for any $\kappa >0$ we proved the existence of a solution $\vv^\kappa$ time-periodic with period $T$ corresponding to initial data  $\vv^\kappa_0$ whose $L^2$-norm is bounded by the constant $K$ given in \eqref{assumption1}. Since we may assume $\kappa<1$,
$\|\vv^\kappa_0\|_{L^2(\O; R^3)}$ is uniformly bounded in $\kappa$, hence there exists a subsequence, still denoted by $\vv^\kappa_0$, weakly converging in $L^2(\O; R^3)$ to a limit $\vv_0$ such that 
\begin{equation}\label{K-bar}
\|\vv_0\|_{L^2(\O; R^3)}\leq \overline{K}:= \left( \frac{C_2(\max_{[0,T]}\|\bfb\|_{L^2(\Omega_T; R^3)}^{q'}+ 1)}{C_1C_S}\right)^{\frac{1}{q}}.
\end{equation}
Next, the sequence $\{\vv^\kappa\}$ satisfies the bounds in \eqref{uniform1} uniformly in $\kappa$. Since
\begin{equation*}\begin{split}
&\intTO{((\kappa+|\Ds\vv|^2)^{\frac{q-2}{2}} |\Ds\vv|)^{q'}}\\
&\leq \intTO{(\kappa+|\Ds\vv|^2)^{\frac{(q-2)q}{2(q-1)}}(k+|\Ds\vv|^2)^{\frac{q}{2(q-1)}}} \\&
= \intTO{(\kappa+|\Ds\vv|^2)^{\frac{q}{2}}}\\&\leq 
\int_{|\Ds\vv|^2\geq \kappa} (2|\Ds\vv|^2)^{\frac{q}{2}}\dx\dt + \int_{|\Ds\vv|^2\leq \kappa} (2\kappa)^{\frac{q}{2}}\dx\dt, 
\end{split}\end{equation*}
we can state that $\S^\kappa:=(\kappa+|\Ds\vv^\kappa|^2)^{\frac{q-2}{2}} \Ds\vv^\kappa$ is uniformly bounded in $L^{q'}(\O_T; R^{3\times 3})$, moreover
 $\vv^\kappa$ enjoys also the bounds in \eqref{uniform2} and the sequence 
$\{\partial_t\vv^\kappa\}$ is bounded in $(L^{\frac{11}{5}}(0,T; W_0^{1, \frac{11}{5}}(\Omega; R^3)))^*$ uniformly in $\kappa$. Hence, employing also the Aubin-Lions compactness lemma, there exist subsequences $\vv^\kappa$, $\S^\kappa$, not relabelled, weakly, weakly-* and strongly converging in the relevant spaces to the limit $\vv$ and $\S$ respectively. We may follow the steps of the passage to the limit in the previous paragraph and obtain that 
\begin{equation}\label{weak-ep}\begin{split}
&\int_0^T \langle\partial_t \vv, \bph\rangle \dt +\intTO{\S:\Ds\bph} - \intTO{\vv\otimes\vv:\Ds\bph} \\
&+\varepsilon \intTO{\Grad \vv:\Grad\bph} + \ep \intTO{|\Ds\vv|^{\frac{1}{5}} \Ds\vv:\Ds\bph}
\\&= \intTO{\bfb \cdot\bph} \mbox{ for any } \bph \in L^{\frac{11}{5}}(0,T; W_{0, {\rm div}}^{1, \frac{11}{5}}(\Omega; R^3)).
 \end{split}\end{equation}
It remains to be proved the identification of the limit $\S$. 
Since 
$$|(\kappa+|\Ds\bph|^2)^{\frac{q-2}{2}}\Ds\bph- |\Ds\bph|^{q-2}\Ds\bph|^{q'}\to 0 \mbox{ as } \kappa\to 0, \mbox{a.e. in } \O_T, $$
and 
$$|(\kappa+|\Ds\bph|^2)^{\frac{q-2}{2}}\Ds\bph- |\Ds\bph|^{q-2}\Ds\bph|^{q'}\!\!\leq \!\! (C+2) |\Ds\bph|^{(q-2)q'}|\Ds\bph|^{q'}\!\leq \!C|\Ds\bph|^q $$
where $C$ is a constant independent of $\kappa$ and varying from line to line, the Lebesgue dominated convergence theorem ensures that 
\begin{equation}\label{lebesgue}
(\kappa+|\Ds\bph|^2)^{\frac{q-2}{2}}\Ds\bph\to |\Ds\bph|^{q-2}\Ds\bph \mbox{ strongly in } L^{q'}(\O_T; R^{3\times 3}).
\end{equation} 
The monotonicity of the operator \eqref{eq-S} implies that the quantity
$$\intTO{[(\kappa+|\Ds\vv^\kappa|^2)^{\frac{q-2}{2}} \Ds\vv^\kappa - (\kappa+|\Ds\bph|^2)^{\frac{q-2}{2}}\Ds\bph]: (\Ds\vv^\kappa-\Ds\bph)\!\!}$$
is non-negative and, using \eqref{system}, it can be rewritten as 
\begin{equation}\label{S-kappa}\begin{split}
&\!\!\!\!\!\int_0^T\!\!\!\intO{\left[\bfb\cdot \vv^\kappa -\varepsilon|\Grad \vv^\kappa|^2-\ep |\Ds\vv^\kappa|^{\frac{11}{5}} -(\kappa+|\Ds\vv^\kappa|^2)^{\frac{q-2}{2}} \Ds\vv^\kappa:\Ds\bph \right.\\
&\left.- (\kappa+|\Ds\bph|^2)^{\frac{q-2}{2}}\Ds\bph: (\Ds\vv^\kappa-\Ds\bph)\right]\!\!}\!\dt.
\end{split}\end{equation}
Employing the strong convergence \eqref{lebesgue}, the boundedness of $\|\Ds\vv^\kappa\|_{L^q(\O_T; R^{3\times 3})}$ uniformly in $\kappa$, and the weak convergence of $\Ds\vv^\kappa$ to $\Ds\vv$ in $L^q(\O_T; R^{3\times 3})$, it follows that 
\begin{equation}\begin{split}
&\intTO{(\kappa+|\Ds\bph|^2)^{\frac{q-2}{2}}\Ds\bph: (\Ds\vv^\kappa-\Ds\bph)} \\
&= \intTO{\left((\kappa+|\Ds\bph|^2)^{\frac{q-2}{2}}\Ds\bph - |\Ds\bph|^{q-2}\Ds\bph\right): \Ds\vv^\kappa}\\
&+ \intTO{|\Ds\bph|^{q-2}\Ds\bph : \Ds\vv^\kappa} \\
&- \intTO{(\kappa+|\Ds\bph|^2)^{\frac{q-2}{2}}\Ds\bph :\Ds\bph}\\
&\longrightarrow  \intTO{|\Ds\bph|^{q-2}\Ds\bph :(\Ds\vv-\Ds\bph)}. 
\end{split}\end{equation}
Then taking the liminf in \eqref{S-kappa} and using all the collected convergences we get
\begin{equation}\label{anna}\begin{array}{l}\displaystyle\vspace{6pt}
0\leq \!\!\int_0^T\!\!\intO{\left[\bfb\cdot \vv -\varepsilon|\Grad \vv|^2-\ep |\Ds\vv|^{\frac{11}{5}}\right.
\\ \displaystyle \vspace{6pt}
\left. -\S:\Ds\bph -|\Ds\bph|^{q-2}\Ds\bph \!:\!(\Ds\vv-\Ds\bph)\right]\!\!}\!\dt.
\end{array}\end{equation}
At this level of approximation, due to the presence of $\ep>0$, we can still use the weak solution itself $\vv$ as test function in \eqref{weak-ep} and comparing the outcome with \eqref{anna}, it follows that 
\begin{equation}
0\leq \intTO{(\S- |\Ds\bph|^{q-2}\Ds\bph) : (\Ds\vv-\Ds\bph)}.
\end{equation}
Now, choosing $\bph:= \vv\pm \lambda {\bf w}$ with $\lambda>0$, dividing by $\lambda$ and then taking the limit as $\lambda\to 0$ we obtain that 
\begin{equation}\label{S-ep}
\S= |\Ds\vv|^{q-2}\Ds\vv \mbox{ a.e. in } \O_T.
\end{equation}

\paragraph{Limit as $\ep\to0$.}
We constructed a sequence $\{\vv^\varepsilon\}$ of time-periodic solution with period $T$ to \eqref{weak-ep} with viscous stress tensor $\{\S^\ep\}$ fulfilling \eqref{S-ep}, corresponding to a sequence of initial data $\{\vv_0^\ep\}$ satisfying \eqref{K-bar}. Taking as test function $\bph=\vv^\varepsilon$ in \eqref{weak-ep} we derive the following uniform estimates
\begin{equation}\label{uniform65-0}\begin{split}
\sup_{\ep>0}&\left(\|\Ds\vv^\varepsilon\|_{L^q(\Omega_T; R^{3\times 3})}+\|\varepsilon^{\frac{1}{2}}\nabla\vv^\varepsilon\|_{L^2(\Omega_T; R^{3\times 3})}\right.	\\
&	\left.+\|\varepsilon^{\frac{5}{11}}\Ds\vv^\varepsilon\|_{L^{\frac{11}{5}}(\Omega_T; R^{3\times 3}))}\right)<+\infty 
\end{split}\end{equation}
As a consequence and similarly to the Section $q\geq 11/5$ we get 
\begin{equation}\label{uniform65}
\sup_{\ep>0}\left(\|\S^\ep\|_{L^{q'}(\Omega_T; R^{3\times3})}+ \|\vv^\varepsilon\|_{L^{\frac{5q}{3}}(\Omega_T; R^3)}\right)<+\infty
\end{equation}
and then
\begin{equation}
\sup_{\ep>0}\|\partial_t\vv^\varepsilon\|_{(L^{\frac{5q}{5q-6}}(0,T; W_0^{1, \frac{5q}{5q-6}}(\Omega; R^3)))^*}<+\infty. 
\end{equation}
Repeating the same arguments perfomed in the Section $q\geq 11/5$, we obtain the existence of $\vv$ and $\S$ such that, for suitable subsequences not relabelled, the following convergences hold:
\begin{align}
\vv^\varepsilon &\to \vv &&\mbox{ weakly-* in } L^\infty(0,T; L^2(\Omega; R^3)),\label{convergence1}\\
 \vv^\varepsilon &\to \vv &&\mbox{ weakly  in } L^q(0,T; W^{1, q}_{0, {\rm div}}(\Omega; R^3)),\\
\partial_t\vv^\varepsilon&\to \partial_t\vv &&\mbox{ weakly-* in } (L^{\frac{5q}{5q-6}}(0,T; W_{0, {\rm div}}^{1, \frac{5q}{5q-6}}(\Omega; R^3)))^*,\\
\vv^\varepsilon&\to \vv &&\mbox{ strongly in } L^r(\Omega_T; R^3) \mbox{ for any } r< \frac{5q}{3},\label{strong}\\
\S^\varepsilon&\to\bar{\S} &&\mbox{ weakly  in } L^{q'}(\Omega_T; R^{3\times 3}).\label{convergence5}
\end{align}
Furthermore,
\begin{equation}\begin{split}
\left|\varepsilon\!\!\int_0^T\int_\Omega\Grad\vv^\varepsilon:\Grad\bph \dx\dt\right|&\leq \varepsilon\, \|\Grad\vv^\varepsilon\|_{\frac{5q}{6}}\|\Grad\bph\|_{\frac{5q}{5q-6}} \\
&\leq \ep \,C\, \|\Grad\bph\|_{\frac{5q}{5q-6}}\to 0 \mbox{ as } \varepsilon \to 0,\end{split}
\end{equation}
and 
\begin{equation}\label{ep-2}\begin{split}
&\!\!\!\!\!\left|\varepsilon\!\int_0^T\!\!\!\int_\Omega|\Ds\vv^\varepsilon|^{\frac{1}{5}}\Ds\vv^\varepsilon:\Ds\bph\dx\dt \right| \leq \,\ep^{\frac{5}{11}}\!\!\int_0^T\!\!\!\intO{\ep^{\frac{6}{11}}|\Ds\vv^\ep|^{\frac{6}{5}}|\Ds\bph|}\dt\\
&\leq \,\ep^{\frac{5}{11}} \left(\intTO{\ep |\Ds\vv^\ep|^{\frac{11}{5}}}\right)^{\frac{6}{11}} \|\Ds\bph\|_{L^{\frac{11}{5}}(\O_T; R^{3\times 3})}\\&= \varepsilon^{\frac{5}{11}} \|\varepsilon^{\frac{5}{11}}\Ds\vv^\varepsilon\|_{\frac{11}{5}}^{\frac{6}{5}}\|\Ds\bph\|_{\frac{11}{5}} \leq \varepsilon^{\frac{5}{11}} \, C\,\|\Ds\bph\|_{\frac{11}{5}} \to 0 \mbox{ as } \varepsilon \to 0.
\end{split}\end{equation}
Then taking the limit as $\ep\to 0$ in \eqref{weak-ep} and by a standard density argument we obtain that $\vv$ and $\bar{\S}$ fulfil the following weak formulation:
\begin{equation}\label{final-weak}\begin{split}
\int_0^T \langle\partial_t \vv, \bph\rangle \dt - \intTO{\vv\otimes\vv: \Grad \bph} + \intTO{\bar{\S}:\Ds\bph} &\\= \intTO{\bfb \cdot\bph}
\mbox{ for any } \bph\in C^\infty((0, T); C_{0, {\rm div}}^\infty(\O; R^3)).&
 \end{split}\end{equation}
Being $(\vv^\varepsilon(t), \bph)$ with $\bph\in C_{0, {\rm div}}^\infty(\O; R^3)$ uniformly bounded and equicontinuous thanks to \eqref{uniform65-0} and \eqref{uniform65}, the Ascoli-Arzel\`{a} theorem gives the following convergence for a suitable subsequence:
\begin{equation}
(\vv^\varepsilon(t), \bph)\to (\vv(t), \bph) \mbox{ uniformly in } t\in[0,T],  \mbox{ as } \varepsilon\to 0
\end{equation}
and, $(\vv(t), \bph)$ is a continuous function thus
$\vv\in C_{\rm weak}(0,T; L^2(\Omega; R^3)).$ 
On the other hand, by \eqref{K-bar} we get that  
\begin{equation}\label{weak}
\vv_0^\varepsilon\to \vv_0 \mbox{ weakly  in } L^2(\Omega; R^3)
\end{equation}
for a subsequence not relabelled. Therefore it holds 
\begin{equation}
(\vv(T)-\vv_0, \bph)=\lim_{\varepsilon\to 0} (\vv^\varepsilon(T)-\vv_0^\varepsilon, \bph)=0 \ \ \forall \bph\in C^\infty_{0, {\rm div}}(\Omega; R^3),
\end{equation}
which means that the limit function $\vv$ is time-periodic with period $T$. It remains to identify the rheological law. We perform the Lipschitz truncation method in the solenoidal version. By \eqref{final-weak}, from the weak formulation fulfilled by the approximation $\vv^\ep$ i.e. \eqref{weak-ep}, and taking into account the time-periodicity of $\vv^\varepsilon$ and $\vv$, it is easily seen that 
\begin{equation}\label{system2}\begin{array}{l}\displaystyle\vspace{6pt}
-\int_0^T\int_\Omega (\vv^\varepsilon-\vv)\cdot\partial_t\bph\dx\dt+\int_0^T\int_\Omega (\S^\varepsilon-\bar{\S}):\Ds\bph\dx\dt\\\displaystyle\vspace{6pt}
+\varepsilon\int_0^T\int_\Omega\Grad \vv^\varepsilon: \Grad\bph \dx\dt
+\, \varepsilon\int_0^T\int_\Omega|\Ds\vv^\varepsilon|^{\frac{1}{5}}\Ds\vv^\varepsilon:\Ds\bph\dx\dt\\\displaystyle\vspace{6pt}
 =\int_0^T\int_\Omega(\vv^\varepsilon\otimes\vv^\varepsilon-\vv\otimes\vv):\Grad\bph \dx\dt
\end{array}
\end{equation}
 for any $\bph\!\in\! C^\infty((0,T); C^\infty_{0, {\rm div }}(\Omega; R^3))$. Let us set $\vu^\varepsilon:=\vv^\varepsilon-\vv$ and $\mathbb{H}^\varepsilon_1:=\S^\varepsilon-\bar{\S}$, then let us rewrite convergences \eqref{convergence1}--\eqref{convergence5} as follows:
\begin{align}
\vu^\varepsilon &\to 0 &&\mbox{ weakly-* in } L^\infty(0,T; L^2(\Omega; R^3)),\\
\vu^\varepsilon &\to 0 &&\mbox{ weakly  in } L^q(0,T; W^{1, q}_{0, {\rm div }}(\Omega; R^3)),\\
\vu^\varepsilon &\to 0 &&\mbox{ strongly in } L^r(\Omega_T; R^3)) \mbox{ for any } 1<r< \frac{5q}{3},\\
\mathbb{H}^\varepsilon_1&\to 0 &&\mbox{ weakly  in } L^{q'}(\Omega_T; R^{3\times 3}).
\end{align}
Further, since for any $1<r<5q/6$ it holds 
\begin{align*}
\int_0^T\int_\Omega |\varepsilon\Grad \vv^\varepsilon|^r\dx\dt&=\int_0^T\int_\Omega \varepsilon^{\frac{r}{2}+ \frac{r}{2}}|\Grad\vv^\varepsilon|^r \dx\dt \\ &
\leq  \varepsilon^{\frac{r}{2}} \left(\int_0^T\int_\Omega \varepsilon|\Grad\vv^\varepsilon|^2 \dx\dt\right)^{\frac{2}{r}}|\Omega_T|^{\frac{2}{2-r}},\\
\int_0^T\int_\Omega (\varepsilon|\Ds \vv^\varepsilon|^\frac{6}{5})^r \dx\dt&=\int_0^T\int_\Omega \varepsilon^{\frac{5r}{11}+ \frac{6r}{11}}|\Ds\vv^\varepsilon|^\frac{6r}{5} \dx\dt\\&
 \leq  \varepsilon^\frac{5r}{11} \left(\int_0^T\int_\Omega \varepsilon |\Ds\vv^\varepsilon|^\frac{11}{5} \dx\dt\right)^{\frac{6r}{11}}|\Omega_T|^{\frac{11}{11-6r}},
\end{align*}
thus using also \eqref{uniform1} we get that for any $ r< {5q}/{6}$
\begin{equation}\label{strong-ep} \begin{split}
\mathbb{H}^\varepsilon_2&:=\varepsilon\Grad \vv^\varepsilon+\varepsilon|\Ds\vv^\varepsilon|^{\frac{1}{5}}\Ds\vv^\varepsilon +(\vv\otimes\vv-\vv^\varepsilon\otimes\vv^\varepsilon)\\
& \to 0
  \mbox{ strongly in } L^r(\Omega_T; R^{3\times 3}) \mbox{ as } \varepsilon \to 0.
\end{split}\end{equation}
Therefore, \eqref{system2} can be rewritten again as:
\begin{equation}
\int_0^T\int_\Omega \vu^\varepsilon\cdot\partial_t\bph\dx\dt=\int_0^T\int_\Omega (\mathbb{H}^\varepsilon_1 +\mathbb{H}^\varepsilon_2) :\Ds\bph\dx\dt
\end{equation}
for any  $\bph\in C^\infty(0,T; C^\infty_{0, {\rm div }}(\Omega; R^3))$. At this point let us set $\ep:=\ep_m=1/m$. Then all the assumptions of Theorem \ref{lipschitz-trunc} and Corollary \ref{lipschitz-trunc2} in Section 4 are fulfilled. Particularly, note that we can consider any small $r>1$. Now, take any cylinder $Q\subset (0, T)\times \Omega$, for suitable $\xi\in C_0^\infty(\frac{1}{6}Q; R^3)$,\footnote{For any $\alpha>0$ we denote $\alpha Q$ the cylinder $Q$ scaled by $\alpha$ with respect to its center.} for any $\theta\in (0, 1)$, set $\S:= |\Ds\vv|^{q-2}\Ds\vv$, then consider the following quantity 
$$\begin{array}{c}\displaystyle \vspace{6pt}
\int_Q((\S^m-\S): (\Ds\vv^m-\Ds\vv))^\theta \xi \dx\dt \\ \displaystyle\vspace{6pt}
= \int_Q((\S^m-\S): (\Ds\vv^m-\Ds\vv))^\theta \xi \,\chi_{O_{m,k}}\dx\dt  \\ \displaystyle+ \int_Q((\S^m-\S): (\Ds\vv^m-\Ds\vv))^\theta \xi \,\chi_{O_{m,k}^C}\dx\dt.
\end{array}$$
Let us treat the two terms separately. For the first term using the H\"{o}lder inequality, the established uniform estimates in $\ep_m$, the statements (a) and (h) of Theorem \ref{lipschitz-trunc} we obtain
\begin{equation}\label{lip1}\begin{array}{l}\displaystyle \vspace{6pt}
\limsup_{m\to +\infty} \int_Q((\S^m-\S): (\Ds\vv^m-\Ds\vv))^\theta \xi \,\chi_{O_{m,k}}\dx\dt \\ \displaystyle\vspace{6pt}
\leq C\limsup_{m\to +\infty}\!\left(\int_Q\!(\S^m-\S)\!:\! (\Ds\vv^m-\Ds\vv) \dx\dt\right)^{\theta}\!\! |O_{m,k}|^{1-\theta}
\\ \displaystyle\leq C 2^{k(\theta-1)}.
\end{array}\end{equation}
On the other hand, in the second term H\"{o}lder's inequality gives
$$\begin{array}{c}\displaystyle \vspace{6pt}
\int_Q((\S^m-\S): (\Ds\vv^m-\Ds\vv))^{\theta} \xi \,\chi_{O_{m,k}^C}\dx\dt \\ \displaystyle \vspace{6pt}\hfill =\int_{O_{m,k}^C}((\S^m-\S): (\Ds\vv^m-\Ds\vv)\,\xi)^{\theta} \,\xi^{1-\theta}\dx\dt\\\displaystyle \vspace{6pt}\hfill
\leq C\left(\int_Q(\S^m-\S): (\Ds\vv^m-\Ds\vv)\,\xi \,\chi_{O_{m,k}^C}\dx\dt  \right)^{\theta}. \end{array}$$
Setting $\mathbb{K}:= \bar{\S}-\S$ it results true the identity 
$$\S^m-\S=\mathbb{H}^m_1+\mathbb{K},$$
and employing Corollary \ref{lipschitz-trunc2} it follows that 
\begin{equation}\label{lip2}
\limsup_{m\to +\infty}\int_Q((\S^m-\S): (\Ds\vv^m-\Ds\vv))^{\theta}\, \xi \,\chi_{O_{m,k}^C}\dx\dt   \leq C\,2^{-\frac{k\theta}{q}}.
\end{equation}
Gathering \eqref{lip1} and \eqref{lip2} it results that
$$\limsup_{m\to +\infty}\int_Q((\S^m-\S): (\Ds\vv^m-\Ds\vv))^\theta\, \xi \dx\dt\leq C\, 2^{k(\theta-1)}+ C\,2^{-\frac{k\theta}{q}},$$
and since $k$ can be chosen arbitrarily large  finally we have that 
\begin{equation}\begin{array}{c}\displaystyle \vspace{6pt}
\limsup_{m\to +\infty}\int_{\frac{1}{8}Q} \left((\S^m-\S): (\Ds\vv^m-\Ds\vv)\right)^\theta \dx\dt\\\displaystyle \vspace{6pt}
\leq \limsup_{m\to +\infty}\int_Q((\S^m-\S): (\Ds\vv^m-\Ds\vv))^\theta\, \xi \dx\dt=0.
\end{array}\end{equation}
By the monotonicity of the integrand on the left-hand side it follows that 
\begin{equation}
\lim_{m\to+\infty} (\S^m-\S): (\Ds\vv^m-\Ds\vv) =0 \mbox{ a.e. in } \frac{1}{8}Q, 
\end{equation}
thus by virtue of the results of Dal Maso and Murat \cite{DalMasoMurat} we conclude that
\begin{equation}
\lim_{m\to+\infty} \Ds\vv^m= \Ds\vv  \mbox{ a.e. in } \frac{1}{8}Q,
\end{equation}
and then the continuity of the stress tensor implies that 
\begin{equation}
\lim_{m\to+\infty} \S^m=\S= |\Ds\vv|^{q-2}\Ds\vv \mbox{ a.e. in } \frac{1}{8}Q.
\end{equation}
Since $Q$ was arbitrary we can conclude that 
\begin{equation}
\lim_{m\to+\infty} \S^m=\S= |\Ds\vv|^{q-2}\Ds\vv \mbox{ a.e. in } \Omega_T.
\end{equation}
Being $\S^m$ uniformly integrable, the Vitali convergence theorem yields $\S^m\to \S$ weakly in $L^{q'}(\O_T; R^{3\times 3})$ and finally
$$ \bar{\S}=\S= |\Ds\vv|^{q-2}\Ds\vv \mbox{ a.e. in } \Omega_T.$$
The proof is complete.

\section{Proof of Corollary \ref{corollary}}
By virtue of the proof of Theorem \ref{thm}, we can claim that for any $\varepsilon >0$ there exists $\vv^\varepsilon$ and $\S^\varepsilon$ such that 
\begin{equation}\label{cor1}\begin{split}
\!\!\!\!\int_0^\tau \langle\partial_t \vv^\varepsilon, \bph\rangle \dt - \int_0^\tau\intO{\vv^\varepsilon\otimes\vv^\varepsilon: \Grad \bph}\dt + \int_0^\tau\intO{\S^\varepsilon:\Ds\bph}\dt &\\
+\varepsilon\int_0^\tau \intO{\Grad \vv^\varepsilon:\Grad\bph}\dt + \ep \int_0^\tau\intO{|\Ds\vv^\varepsilon|^{\frac{1}{5}} \Ds\vv^\varepsilon:\Ds\bph}\dt&
\\
= \int_0^\tau\intO{\bfb \cdot\bph}\dt
\mbox{ for any } \bph\in L^{\frac{11}{5}}(0,T; W_{0, {\rm div}}^{1, \frac{11}{5}}(\Omega; R^3)),&
 \end{split}\end{equation}
for any $0<\tau\leq T$, with 
$$ \S^\varepsilon=|\Ds\vv^\varepsilon|^{q-2}\Ds\vv^\varepsilon \mbox{ a.e. in } \Omega_T,$$
and
$$ \vv^\varepsilon(T)=\vv^\varepsilon(0) \mbox{ in } L^2(\Omega; R^3).$$
Let us consider $\bph=\vv^\varepsilon$ in \eqref{cor1} it follows that  
\begin{equation}\begin{split}
\frac{1}{2}\Dt \|\vv^\varepsilon(t)\|_{L^2(\Omega; R^3)}^2 +& \|\Ds\vv^\varepsilon(t)\|_{L^q(\Omega; R^{3\times 3})}^q\\
&\leq \|\bfb(t)\|_{L^2(\Omega; R^3)}\|\vv^\varepsilon(t)\|_{L^2(\Omega; R^3)}\\
&\leq C\|\bfb(t)\|_{L^2(\Omega; R^3)}^{q'} + \frac{1}{2}\|\Ds\vv^\varepsilon(t)\|_{L^q(\Omega; R^{3\times 3})}^q
\end{split}\end{equation}
for a.a. $t\in (0, T)$, then 
\begin{equation}\label{cor2}
\frac{1}{2}\Dt \|\vv^\varepsilon(t)\|_{L^2(\Omega; R^3)}^2 + \alpha \|\vv^\varepsilon(t)\|_{L^2(\Omega; R^3)}^q \leq C\|\bfb(t)\|_{L^2(\Omega; R^3)}^{q'}
\end{equation}
for a.a. $t\in (0, T)$, where $\alpha >0$ is the constant due to the Sobolev embedding $W_0^{1, q}\hookrightarrow L^2$, while $C>0$ is a constant that may vary from line to line but it is independent of $\varepsilon$ and of $t$. Since $\|\bfb(t)\|_{L^2(\Omega; R^3)}=0$ in $[\bar{t}, T]$, \eqref{cor2} implies that 
\begin{equation}
\frac{1}{2}\Dt \|\vv^\varepsilon(t)\|_{L^2(\Omega; R^3)}^2 + \alpha \|\vv^\varepsilon(t)\|_{L^2(\Omega; R^3)}^q \leq 0 \mbox{ for a.a. } t\in (\bar{t}, T).
\end{equation}
Restricting $q\in (6/5; 2)$, integrating between $\bar{t}$ and $t$, and employing \eqref{K-bar} we get 
\begin{equation}\begin{split}
 \|\vv^\varepsilon(t)\|_{L^2(\Omega; R^3)}^{2-q} &\leq  \|\vv^\varepsilon(0)\|_{L^2(\Omega; R^3)}^{2-q} - (2-q) \alpha(t-\bar{t})\\
& \leq \overline{K}^{2-q} - (2-q) \alpha(t-\bar{t}),
\end{split} \end{equation}
therefore
\begin{equation}
 \|\vv^\varepsilon(t)\|_{L^2(\Omega; R^3)}=0 \mbox{ for a.a. } t\mbox{: } \bar{t} + \frac{\overline{K}^{2-q}}{\alpha (2-q)} \leq t\leq T.
 \end{equation}
Now let us follow the limit as $\varepsilon \to 0$ as showed in the previous section. In particular it holds
$$\vv^\varepsilon(t)\to \vv(t) \mbox{ weakly in } L^2(\Omega; R^3) \mbox{ for a.a. } t\in (\bar{t}_{v}, T)$$
where $$\bar{t}_v:=\bar{t} + \frac{ \overline{K}^{2-q}}{\alpha (2-q)},$$ 
then the weak-lower semicontinuity of the $L^2$-norm gives that 
$$\|\vv(t)\|_{L^2(\Omega; R^3)}\leq \liminf_{\varepsilon \to 0 } \|\vv^\varepsilon(t)\|_{L^2(\Omega; R^3)} = 0 \mbox{ for a.a. } t\in (\bar{t}_v, T)$$
and this finishes the proof.

\section{Auxilliary tools}
 We state the divergence-free Lipschitz truncations of Bochner-Sobolev functions taken from \cite{BreDieSch}.
 \begin{Theorem}\label{lipschitz-trunc}
 Let $Q=I\times B\subset R\times R^3$ be a space-time cylinder and let $1<q<+\infty$ with $q, q':=q/(q-1) >r>1$. Let $\{\vu^m\}$ and $\{\mathbb{H}^m\}$ with $\mathbb{H}^m:=\mathbb{H}^m_1+\mathbb{H}^m_2$ be sequences fulfilling
 \begin{align}
& \Div\vu^m =0 \mbox{ a.e. in } Q,\\
& \partial_t \vu^m=-\Div \mathbb{H}^m \mbox{ in sense of distributions in } C^\infty(I; C_{0, {\rm div}}^\infty(B; R^3)),
 \end{align}
assume that $\vu^m$ is bounded in $L^\infty(I; L^r(B; R^3))$ uniformly and  that 
 \begin{align}
 \vu^m &\to 0 &&\mbox{ weakly  in } L^q(I; W^{1, q}_0(B; R^3)),\\
\vu^m &\to 0 &&\mbox{ strongly in } L^r(Q; R^3)),\\
\mathbb{H}^m_1&\to 0 &&\mbox{ weakly  in } L^{q'}(Q; R^{3\times 3}),\\
\mathbb{H}^m_2&\to 0  &&\mbox{ strongly in } L^r(Q; R^{3\times 3})
 \end{align}
 as $m\to + \infty$. Then, there is a double sequence $\{ \lambda_{m, k}\}_{m, k=1}^\infty\subset (0, \infty)$ such that 
 \begin{itemize}
 \item[(a)] $2^{2^k}\leq \lambda_{m, k}\leq 2^{2^{k+1}},$
 \end{itemize}
 there exist a double sequence of functions $\{\vu^{m,k}\}_{m, k=1}^\infty\subset L^1(Q; R^3)$, a double sequence $\{O_{m,k}\}_{m, k=1}^\infty$ of measurable subsets of $Q$, a constant $C>0$ and $k_0\in N$ such that for any $k\geq k_0$ the following properties are satisfied:
 \begin{itemize}
\item[(b)] $\vu^{m,k}\in L^s(\frac{1}{4} I; W^{1, s}_{0, {\rm div}}(\frac{1}{6}B; R^3))$ for any $s\in (1, \infty)$,
\item[(c)] ${\rm supp}\, \vu^{m,k}\subset \frac{1}{6} Q$,
\item[(d)] $\vu^{m,k}=\vu^m$ a.e. in $\frac{1}{8}Q\setminus O_{m,k}$, 
\item[(e)] $\|\Grad\vu^{m,k}\|_{L^\infty(\frac{1}{4}Q; R^{3\times3})}\leq C\lambda_{m,k}$,
\item[(f)] $\vu^{m,k} \to 0$ strongly in  $L^\infty(\frac{1}{4}Q; R^3)$ as $m\to \infty$,
\item[(g)] $\Grad\vu^{m,k} \to 0$ weakly-* in $L^\infty(\frac{1}{4}Q; R^3)$ as $m\to\infty$,
\item[(h)] $\limsup_{m\to +\infty}(\lambda_{m,k})^q|O_{m,k}|\leq C 2^{-k}$,
\item[(i)] $\limsup_{m\to +\infty}\left|\int_{\frac{1}{8}Q\setminus O^{m,k}} \mathbb{H}^m:\Grad \vu^{m,k} \dx\dt\right|\leq C(\lambda_{m,k})^q|O_{m,k}|$.
\end{itemize}  
\end{Theorem}
 \begin{Corollary}\label{lipschitz-trunc2}
 Let all the assumptions of Theorem \ref{lipschitz-trunc} be fulfilled and assume that $\vu^m$ is bounded in $L^\infty(I; L^r(B; R^3))$ uniformly in $m$. Then there exist $\xi\in C_0^\infty(\frac{1}{6}Q; R^3)$ such that $\chi_{\frac{1}{8}Q}\leq \xi\leq \chi_{\frac{1}{6}Q}$  and a constant $C$ depending only on $\xi$ such that for every $\mathbb{K}\in L^{q'}(\frac{1}{6}Q; R^{3\times 3})$ there holds:
\begin{equation}
 \limsup_{m\to+\infty} \left| \int(\mathbb{H}^m_1 + \mathbb{K}): \Grad\vu^m\xi\chi_{O_{m,k}^c}\dx\dt\right|\leq C\,2^{-\frac{k}{q}},
 \end{equation}
 and 
 \begin{align}
& \!\!\vu^{m,k} \to 0 \mbox{ strongly in } L^s(\frac{1}{4}Q; R^3) \mbox{ as } m\to \infty \mbox{ and } k \mbox{ fixed, } \forall s\!\in\! (1,\infty),\\
&\!\!\Grad\vu^{m,k} \!\to \!0 \mbox{ weakly in } L^s(\frac{1}{4}Q; R^3)  \mbox{ as } m\!\to\!\infty \mbox{ and } k \mbox{ fixed, } \forall s\!\in\! (1,\infty).
 \end{align}
 \end{Corollary}
 \section*{Acknowledgements} 
 The research of A.A. is supported by Einstein Foundation, Berlin. She is also member of the Italian National Group for the Mathematical Physics (GNFM) of INdAM. 

 \section*{Compliance with ethical standards}
 The author states that there is no conflict of interest.


\begin{thebibliography}{10}

\bibitem{AnnaEd}
A.~Abbatiello and E.~Feireisl, \emph{On a class of generalized solutions to
  equations describing incompressible viscous fluids}, Ann. Mat. Pura Appl. (4)
  \textbf{199} (2020), no.~3, 1183--1195.

\bibitem{Abba}
A.~Abbatiello and P.~Maremonti, \emph{Electrorheological fluids: ill posedness of uniqueness backward in time}, Nonlinear Anal.\textbf{170} (2018), 47--69.

\bibitem{AbbatielloMaremonti}
A.~Abbatiello and P.~Maremonti, \emph{Existence of regular time-periodic
  solutions to shear-thinning fluids}, J. Math. Fluid Mech. \textbf{21} (2019),
  no.~2, Art. 29, 14.

\bibitem{Barhoun}
A.~Barhoun and A.~B. Lemlih, \emph{A reproductive property for a class of
  non-{N}ewtonian fluids}, Appl. Anal. \textbf{81} (2002), no.~1, 13--38.

\bibitem{BleMalRaj}
J.~Blechta, J.~M\'{a}lek, and K.~R. Rajagopal, \emph{On the {C}lassification of
  {I}ncompressible {F}luids and a {M}athematical {A}nalysis of the {E}quations
  {T}hat {G}overn {T}heir {M}otion}, SIAM J. Math. Anal. \textbf{52} (2020),
  no.~2, 1232--1289.

\bibitem{BreDieSch}
D.~Breit, L.~Diening, and S.~Schwarzacher, \emph{Solenoidal {L}ipschitz
  truncation for parabolic {PDE}s}, Math. Models Methods Appl. Sci. \textbf{23}
  (2013), no.~14, 2671--2700.

\bibitem{Crispo}
F.~Crispo, \emph{A note on the existence and uniqueness of time-periodic
  electro-rheological flows}, Acta Appl. Math. \textbf{132} (2014), 237--250.

\bibitem{CriGriMar}
F.~Crispo, C.~Grisanti, and P.~Maremonti, \emph{Singular p-laplacian parabolic system in exterior domains: higher regularity of
solutions and related properties of extinction and asymptotic behavior in time}, Ann. Sc. Norm. Super. Pisa Cl. Sci. (5) \textbf{19} (2019), no.~3, 913--949.

\bibitem{DalMasoMurat}
G.~Dal~Maso and F.~Murat, \emph{Almost everywhere convergence of gradients of
  solutions to nonlinear elliptic systems}, Nonlinear Anal. \textbf{31} (1998),
  no.~3-4, 405--412.

\bibitem{DiB}
E.~DiBenedetto, \emph{Degenerate parabolic equations}, Universitext. Springer-Verlag, New York, 1993.

\bibitem{DRW}
L.~Diening, M.~R\r{u}\v{z}i\v{c}ka, and J.~Wolf, \emph{Existence of weak
  solutions for unsteady motions of generalized {N}ewtonian fluids}, Ann. Sc.
  Norm. Super. Pisa Cl. Sci. (5) \textbf{9} (2010), no.~1, 1--46.

\bibitem{Lions}
J.~L. Lions, \emph{Sur certaines \'{e}quations paraboliques non lin\'{e}aires},
  Bull. Soc. Math. France \textbf{93} (1965), 155--175.

\bibitem{Malekbook}
J.~M\'{a}lek, J.~Ne\v{c}as, M.~Rokyta, and M.~R\r{u}\v{z}i\v{c}ka, \emph{Weak
  and measure-valued solutions to evolutionary {PDE}s}, Applied Mathematics and
  Mathematical Computation, vol.~13, Chapman \& Hall, London, 1996.

\bibitem{Prouse63}
G.~Prouse, \emph{Soluzioni periodiche dell'equazione di {N}avier-{S}tokes},
  Atti Accad. Naz. Lincei Rend. Cl. Sci. Fis. Mat. Nat. (8) \textbf{35} (1963),
  443--447.

\end{thebibliography}
\end{document}